\documentclass[11pt]{amsart}

\usepackage{amssymb,latexsym}
\def\cbstart{}
\def\cbend{}

\newcommand{\Ot}[1]{\mathcal{O}_{{#1}}}

\newcommand {\legendre}[2]{\genfrac {(}{)}{1pt}{}{#1}{#2}}

\newcommand{\GFpn}[2]{\mathbb{F}_{{#1}^{#2}}}
\newcommand{\GFq}[1]{\mathbb{F}_{#1}}
\newcommand{\kb}{\mathbf{k}}
\newcommand{\K}{\mathbf{K}}

\newcommand{\OK}{\mathcal{O}_K}
\newcommand{\Q}{\mathbb{Q}}

\newcommand{\myjac}[2]{\left(\frac{#1}{#2}\right)}

\newcommand{\galfone}[1]{\mathfrak{w}_{#1}}

\newtheorem{theorem}{Theorem}[section]
\newtheorem{prop}{Proposition}[section]
\newtheorem{coro}{Corollary}[section]
\newtheorem{lemma}{Lemma}[section]

\newtheorem{conjecture}{Conjecture}[section]

\begin{document}

  \title{Computing the cardinality of CM elliptic curves using torsion points}
  \author{F. Morain}
  \address{Laboratoire d'Informatique de l'\'Ecole
  polytechnique (LIX)\\
  F-91128 Palaiseau Cedex\\
  France}
  \email[F. Morain]{morain@lix.polytechnique.fr}
  \thanks{The author is on leave from the French Department of
Defense, D\'el\'egation G\'en\'erale pour l'Armement.}
  \date{\today}

  \begin{abstract}
Let $\mathcal{E}$ be an elliptic curve having complex multiplication by a
given quadratic order of an imaginary quadratic field $\K$. The field
of definition of $\mathcal{E}$ is the
ring class field $\Omega$ of the order. If the prime $p$ splits completely in
$\Omega$, then we can reduce $\mathcal{E}$ modulo one the factors of
$p$ and get a curve $E$ defined over $\GFq{p}$. The trace of the
Frobenius of $E$ is known up to sign and we need a
fast way to find this sign. For this, we propose to use the action of the
Frobenius on torsion points of small order built with
class invariants {\it \`a la} Weber, in a manner reminiscent of the
Schoof-Elkies-Atkin algorithm for computing the cardinality of a given
elliptic curve modulo $p$. We apply our results to the Elliptic Curve
Primality Proving algorithm (ECPP).
  \end{abstract}
  \maketitle


\section{Introduction}

Let $\K$ be an imaginary quadratic field of discriminant $-D$. For
any integer $t$, let $\Ot{t}$ be the order of conductor $t$ of
$\K$, $\Delta_t = - t^2 D$ its discriminant, and $h_t = h(\Delta_t)$
its class number. We denote by $\Omega_t$ the ring class field modulo
$t$ over $\K$. By class
field theory, the extension $\Omega_t/\K$ can be constructed using the
minimal polynomial of the modular function $j$ over a set of
representatives $\{\mathfrak{i}_1, \mathfrak{i}_2, \ldots,
\mathfrak{i}_{h_t}\}$ of the class group $Cl(\Ot{t})$.
An elliptic curve $\mathcal{E}$ of invariant $j(\mathfrak{i}_r)$ is
defined over $\Omega_t$ and has complex multiplication (CM) by $\Ot{t}$.
We denote by $H_{\Delta_t}[j](X)$ the minimal
polynomial of the $j$'s, namely
$$H_{\Delta_t}[j](X) = \prod_{r=1}^{h_t} (X - j(\mathfrak{i}_r))$$
which is known to have rational integer coefficients.

Let $p$ be a rational prime number which splits completely in $\Omega_t$, or
equivalently which is the norm of an integer of $\Omega_t$ (that is $p
= (U^2+D t^2 V^2)/4$ for 
rational integers $U$ and $V$). Then we can reduce $\mathcal{E}$
modulo a prime divisor $\mathfrak{P}$ of $p$ to get an elliptic curve
$E/\GFq{p}$ having CM by $\Ot{t}$. If $\pi$ denotes the
Frobenius of $E$, then it can be viewed as an element of $\Ot{t}$ of norm
$p$, that is (assuming that $\Delta_t \not\in \{-3, -4\}$):
\begin{equation}\label{eqnorm}
\pi = (\pm U \pm t V\sqrt{-D})/2.
\end{equation}
The cardinality of $E(\GFq{p})$ is the norm of $\pi-1$, or more simply
$p+1\mp U$.

The $j$-invariant of $E/\GFq{p}$ is the reduction of one of the
$j(\mathfrak{i}_r)$'s modulo $p$, that is a root of $H(X) =
H_{\Delta_t}[j](X)$ modulo $p$. Building $E$ is done as follows: find
a root $j$ of $H(X)$, and deduce from that the equation of $E$. When
$j \not\in \{0, 1728\}$, we may take any equation $E(j, c)$:
$$Y^2 = X^3 + a_4(j) c^2 X + a_6(j) c^3$$
where $c$ is any element of $\GFq{p}$ and
\begin{equation}\label{defEj}
a_4(j) = \frac{3 j}{1728-j}, \quad a_6(j) = \frac{2 j}{1728-j}.
\end{equation}
We will note $E(j)$ for $E(j, 1)$. If its cardinality is $p+1-a$, then
a curve $E(j, c)$ has cardinality $p+1-\legendre{c}{p} a$
(where $\legendre{a}{b}$ stands for the Legendre symbol). A curve with
$\legendre{c}{p}=-1$ is a twist of $E(j)$.
The problem is now to compute $\#E(j)$ modulo $p$, or equivalently,
fix the sign of $U$ in equation (\ref{eqnorm}).

In the course of implementing the ECPP algorithm
\cite{AtMo93b} or for cryptographic reasons, it is important to
compute this cardinality rapidly. We could of course try both signs of
$U$ yielding cardinalities $m$, find some
random points $P$ on $E(j)$ and check whether $[m] P = O_E$ on
$E$. This approach is somewhat probabilistic and we prefer
deterministic and possibly faster solutions.

In the case where $D$ is fundamental and prime to $6$, the solution is to use
Stark's approach \cite{Stark96}, together with tricks described in
\cite{Morain98a}. This method is efficient, provided we can afford
some precomputations. Note that in the special case where $h_t=1$,
which includes $j = 0, 1728$, one
already knows the answer (see \cite{AtMo93b,JoMo95,PaVe96} and the references
given therein). For $D=20$, we have the isolated result of
\cite{LeMo97} (see also section \ref{ssct:20} below). Since the first
version\footnote{{http://arxiv.org/ps/math.NT/0210173}} of the present
article, Ishii \cite{Ishii04} has given the answer for $D$ of class
numbers $2$ or $3$ and divisible by 3, 4, or 5.

Our approach consists in computing the action of the
Frobenius of the curve on torsion points of small order, using the
techniques of the SEA algorithm \cite{Schoof95}. These points
are obtained using singular values of functions on $X_0(\ell)$ for small
prime $\ell$. This will give us algorithmic solutions to our
motivating problem when $\legendre{-D}{\ell} \neq -1$.

Section 2 describes properties of the modular equations defining
$X_0(\ell)$ for prime $\ell$ and their relations to complex
multiplication over $\overline{\Q}$. In Section 3, we briefly describe
the necessary results used in the SEA algorithm. Section 4 contains
our main contribution. We treat the special cases $\ell=3$ in Section
5 and $\ell=5$ in Section 6. Section 7 describes the very interesting
case of $\ell=7$ and for the sake of completeness that of
$\ell=11$. Section 8 is devoted to the particular case $\ell=2$. We
provide numerical examples for each case. We conclude with remarks on the
use of our results in our implementation of ECPP.

The books \cite{Cox89,Silverman94} are a good introduction to all the
material described above.

\section{Modular curves and class invariants}

\cbstart

\subsection{Modular polynomials}

Let $\ell$ be a prime number. The curve $X_0(\ell)$ parametrizes the
cyclic isogenies of degree $\ell$ associated to an elliptic curve $E$
defined over a field $\kb$. An equation for $X_0(\ell)$ can be
obtained as the minimal polynomial of a modular function $f$
invariant under $\Gamma^0(\ell)$. This {\em modular polynomial}, noted
$\Phi[f](X, J)$ is such that $\Phi[f](f(z), j(z))=0$ for all $z$ such
that $\Im z> 0$, where $j(z)$ is the ordinary modular function.

Dedekind's $\eta$ function is
$$\eta(\tau) =  q^{1/24} \prod_{m \geq 1} (1 - q^m)$$
where $q = \exp(2i\pi\tau)$. It is used to build suitable functions
for $\Gamma^0(\ell)$ (see for instance \cite{Newman57,Newman59}). For
example, if
$$\galfone{\ell}(z) = \frac{\eta(z/\ell)}{\eta(z)}$$
and $s = 12/\gcd(12, \ell-1)$, then $\galfone{\ell}^{2s}$ is a modular
function for $\Gamma^0(\ell)$.
The equations for small prime values of $\ell$ are given in Table
\ref{modeq} (see for instance \cite{Morain95a}).

\begin{table}[hbt]
$$\begin{array}{|r|l|}\hline
 \ell & \Phi[\galfone{\ell}] \\ \hline
 2 & { {\left (X+16\right )^{3}} - J {X}} \\
 3 & { {\left (X+27\right )\left (X+3\right )^{3}} - J {X}} \\
 5 & { {\left ({X}^{2}+10\,X+5\right )^{3}} - J {X}} \\
 7 & { {\left ({X}^{2}+13\,X+49\right )\left({X}^{2}+5\,X+1\right )^{
     3}} - J {X}} \\
\hline
\end{array}$$
\caption{Table of modular equations $\Phi[\galfone{\ell}](X, J)$.\label{modeq}}
\end{table}

Among other classes of functions for other modular groups, we find the
classical functions of Weber:
$$\gamma_2(z) = \sqrt[3]{j(z)}, \quad \gamma_3(z) =
\sqrt{j(z)-1728}$$
for which the corresponding modular equations are quite simple.

\subsection{CM theory}

View the class group $Cl(\Delta_t)$ as a set of reduced quadratic
primitive binary forms of discriminant $\Delta_t$, say $Cl(\Delta_t) =
\{ (A, B, C), B^2 - 4 A C = \Delta_t\}$ with $h_t$ forms in it. For a
given $Q = (A, B, C)$, let $\tau_Q = (-B+\sqrt{\Delta_t})/(2A)$. Then
$j(\tau_Q)$ is an algebraic integer that generates
$\Omega_t/\K$. Moreover, the associated curve $E_Q$ of invariant
$j(\tau_Q)$ has CM by $\Ot{t}$.

Suppose $j(\tau)\in \Omega_t$. If $u$ is some function on some
$\Gamma^0(\ell)$, then the roots of $\Phi[u](X,
j(\tau))$ are algebraic integers. They generate an extension of
$\Omega_t$ of degree dividing $\ell+1$. The striking phenomenon, known for
a long time, is that sometimes these roots lie in $\Omega_t$ itself.
We will note $H_{\Delta_t}[u](X)$ for the minimal polynomial of
the invariant $u$.

Among the simplest results in this direction, we have the following,
dating back to Weber \cite{Weber02}.
Suppose $\alpha$ is a quadratic integer with minimal polynomial
$$A \alpha^2 + B \alpha+ C = 0$$
such that $\gcd(A, B, C)=1$ and $B^2-4 A C = \Delta_t$.
\begin{theorem}
If $3\nmid A$, $3\mid B$, then
$$\Q(\gamma_2(\alpha)) = \left\{\begin{array}{ll}
\Q(j(\alpha)) & \text{ if } 3 \nmid \Delta_t,\\
\Q(j(3 \alpha)) & \text{ if } 3 \mid \Delta_t.\\
\end{array}\right.$$
\end{theorem}
A companion result is:
\begin{theorem}\label{thm:g3}
Suppose $2\nmid A$. We assume that
$$B \equiv \left\{\begin{array}{ll}
0\bmod 4 & \text{ if } 2 \mid \Delta_t,\\
1\bmod 4 & \text{ if } 2 \nmid \Delta_t.\\
\end{array}\right.$$
Then
$$\Q(\sqrt{-D} \gamma_3(\alpha)) = \Q(j(\alpha)), \quad \text{ if } 2
\nmid \Delta_t,$$
$$\Q(\gamma_3(\alpha)) = \Q(j(2 \alpha)), \quad \text{ if } 2 \mid
\Delta_t.$$
\end{theorem}
Finding a complete system of conjugate values for $\gamma_2(\alpha)$
(resp. $\gamma_3(\alpha)$), as well as for a lot of such functions, is
explained in \cite{Schertz02}.

\section{The foundations of the SEA algorithm}

 \subsection{Division polynomials and their properties}

For an elliptic $E$, we let $E[n]$ denote the group of $n$-torsion
points of $E$ (over $\overline{\Q}$). We let $f_n^E(X)$ (or simply
$f_n(X)$) denote the $n$-th division polynomial
whose roots are the abscissae of the $n$-torsion points of $E$.
See \cite{Silverman86} for its definition and
properties. For instance for the
curve $E: Y^2=X^3+a X+b$, the first values are:
$$f_0(X) = 0, f_1(X) = 1, f_2(X) = 1,$$
$$f_3(X) = 3\,{X}^{4}+6\,a\,{X}^{2}+12\,b\,X-a^{2},$$
$$f_4(X) =
2\,{X}^{6}+10\,a\,{X}^{4}+40\,b\,{X}^{3}-10\,a^{2}{X}^{2}-8\,a\,b\,X
-2\,a^{3}-16\,b^{2}.$$
Recurrence relations for computing $f_n$ are given by:
$$f_{2 n} = f_n (f_{n+2} f_{n-1}^2 - f_{n-2} f_{n+1}^2),$$
$$f_{2 n+1} = \left\{\begin{array}{ll}
f_{n+2} f_n^3 - f_{n+1}^3 f_{n-1} (16 (X^3+a X+b)^2) & \text{ if } n
\text{ is odd},\\
&\\
16 (X^3+a X+b)^2 f_{n+2} f_n^3 - f_{n+1}^3 f_{n-1} & \text{ if } n
\text{ is even.}\\
\end{array}\right.$$

Remember that the discriminant of $E$ is $\Delta(E) = -2^4 (4 a^3 + 27
b^2)$. We could not find a reference for the following result, though
it may be classical.
\begin{prop}
Let $m$ be an integer. Then
$$\mathrm{Disc}(f_m) = 
\left\{\begin{array}{ll}
2^4 m^{(m^2-12)/2} (-\Delta)^{(m^2-4)(m^2-6)/24} &
\text{if } m \text{ is even and } \geq 4, \\
(-1)^{(m-1)/2} m^{(m^2-3)/2} (-\Delta)^{(m^2-1)(m^2-3)/24} &
\text{if } m \text{ is odd and }\geq 3. \\
\end{array}\right.$$
\end{prop}
Swan's theorem \cite{Swan62} can be used easily to predict the number
of irreducible factors of $f_n(X)$ other a finite field.

\subsection{Explicit factors of $f_n^E(X)$}

Let $E$ be an elliptic curve.
Suppose that we have some modular polynomial $\Phi[f](X, J)$ for a
function $f$ on $\Gamma^0(\ell)$. Then a root $v$ of $\Phi[f](X,
j(E))$ gives rise to a curve which $\ell$-isogenous to $E$, and
to a factor of $f_{\ell}^E(X)$. This is the essence of
the ideas of Elkies and Atkin that improve Schoof's algorithm
for computing the cardinality of curves over finite fields
\cite{Atkin92b,Schoof95,Elkies98}. The computations can be done
using V\'elu's formulas \cite{Velu71} (see also \cite{Morain95a} for
technicalities related to the actual computations). We end up with a
enable us to compute a factor $g_\ell^{E}(X)$ of $f_\ell^{E}(X)$. 

In the table below, for prime $\ell$, we suppose $v_\ell$ is a root of
$\Phi[\galfone{\ell}](X, j)$ and we give the factor $g_\ell^{E(j)}(X)$ of
$f_\ell^{E(j)}(X)$ that can be obtained in Table \ref{factors}.
\begin{table}[hbt]
$$\begin{array}{|c|l|}\hline
\ell & \text{factor} \\ \hline
2 & \left (v_2-8\right )X+v_2+16,\\ \hline
3 &\left (v_3^{2}+18\,v_3-27\right )X+v_3^{2}+30\,v_3+81,\\ \hline
5 &\left (v_5^{2}+4\,v_5-1\right )^{2}\left (v_5^{2}+
22\,v_5+125\right ){X}^{2} \\
  &+2\,\left (v_5^{2}+4\,v_5-1\right )\left (v_5^{2}+10\,v_5+5\right
  )\left (v_5^{2}+22\,v_5+125\right )X \\
  &+ \left (v_5^{2}+22\,v_5+89\right )\left (v_5^{2}+10
 \,v_5+5\right )^{2},\\
\hline
7 & \left(v_7^{4}+14\,v_7^{3}+63\,v_7^{2}+70\,v_7-7\right )^{3}{X}^{3}\\
  &+3\,\left(v_7^{2}+13\,v_7+49\right )
\left (v_7^{2}+5\,v_7+1\right )
\left(v_7^{4}+14\,v_7^{3}+63\,v_7^{2}+70\,v_7-7\right )^{2}{X}^{2}\\
  &+3\,\left (v_7^{2}+13\,v_7+33\right )
\left (v_7^{2}+13\,v_7+49\right )
\left(v_7^{2}+5\,v_7+1\right )^{2} \\
&\quad\quad\quad \times \left
  (v_7^{4}+14\,v_7^{3}+63\,v_7^{2}+70\,v_7-7\right ) X\\
  &+\left (v_7^{2}+13\,v_7+49\right ) \left (v_7^{2}+5\,v_7+1\right)^3
  (v_7^4+26 v_7^3+219 v_7^2+778 v_7+881)\\
\hline
\end{array}$$
\caption{Factors of $f_\ell^{E(j)}$.\label{factors}}
\end{table}

\subsection{The splitting of $\Phi[f](X, j(E))$ in $\GFq{p}$}

We take the following result from \cite{Atkin92b} (see also
\cite{Schoof95}). Let $\ell$ and $p$ be two distinct primes, and
$E/\GFq{p}$ an elliptic curve. Put $\#E = p+1-U$, $\mathcal{D} =
4p-U^2$. We denote the splitting type of a squarefree polynomial
$P(X)$ by the degrees of its factors. For instance, a polynomial of
degree $4$ having two linear factors and one quadratic factor will be
said to have splitting type $(1)(1)(2)$.
\begin{theorem}\label{modeqsplit}
Let $f$ be a function for $\Gamma_0^{\ell}$ and
put $\Psi(X) \equiv \Phi[f](X, j(E)) \bmod p$.

If $\legendre{-\mathcal{D}}{\ell}=0$, then $\Psi$ splits as $(1)(\ell)$ or
$(1) \cdots (1)$.

If $\legendre{-\mathcal{D}}{\ell}=+1$, then $\Psi$ splits as
$(1)(1)(r) \cdots (r)$ where $r \mid \ell-1$ and $r>1$ if $\ell\neq 2$.

If $\legendre{-\mathcal{D}}{\ell}=-1$, then $\Psi$ splits as $(r) \cdots
(r)$ where $r>1$ and $r \mid \ell+1$.

If $k$ denotes the number of factors of $\Psi$, then $(-1)^k =
\legendre{p}{\ell}$.
\end{theorem}
We can make precise the first part of the theorem as follows:
\begin{theorem}\label{kohel}
Let $p = (U^2+D V^2)/4$. If $\ell \mid V$, then $\Psi(X)$ has $\ell+1$
roots modulo $p$.
\end{theorem}

\noindent
{\em Proof:} See Kohel's thesis \cite{Kohel96}. $\Box$

\subsection{Elkies's ideas}

We briefly summarize Elkies's idea \cite{Elkies98}. Let $\pi$ be the
Frobenius of the curve, sending any point $P = (x, y)$ of
$E(\overline{\GFq{p}})$ to $(x^p, y^p)$.
\begin{theorem}
\label{elkies}
Let $\chi(X)=X^2-UX+p$ denote the characteristic polynomial of the
Frobenius $\pi$ of the elliptic $E$ of cardinality $p+1-U$.
When $\legendre{-\mathcal{D}}{\ell} \neq -1$,
the restriction of $\pi$ to $E[\ell]$ (denoted by $\pi|_{E[\ell]}$) has at
least one eigenvalue. To each eigenvalue $\lambda$ of $\pi|_{E[\ell]}$
corresponds a factor of degree $(\ell-1)/2$ of $f_\ell$. We deduce that
$U \equiv \lambda + p/\lambda \bmod\ell$.
\end{theorem}

We will note $g_{\ell, \lambda}(X)$ the factor of
$f_{\ell}^{E(j)}(X)$ associated to the eigenvalue $\lambda$. Let
$\omega$ denote the order of $\lambda$ modulo $\ell$ and $\sigma =
\omega/2$ if $\omega$ is even and $\omega$ otherwise. With these
notations, one can show the following result:
\begin{prop}
The splitting type of $g_{\ell, \lambda}(X) \bmod p$ is
$(\sigma)(\sigma)\cdots(\sigma)$ with $\kappa$ factors such that
$(\ell-1)/2 = \kappa \sigma$.
\end{prop}
From this, we deduce:
\begin{coro}
The polynomial $g_{\ell, \lambda}(X)$ splits completely modulo $p$ if
and only if $\lambda = \pm 1 \bmod \ell$.
\end{coro}

Note also the following result of Dewaghe \cite{Dewaghe98} in the
formulation of \cite{MaMu01}.
\begin{prop}
Let $r = \mathrm{Resultant}(g_{\ell, \lambda}(X), X^3+a_4(j)X+a_6(j))$. Then
$$\legendre{\lambda}{\ell} = \legendre{r}{p}$$
\end{prop}
Classically, this enables us to fix the sign of $\lambda$ when
$\ell\equiv 3\bmod 4$.

\section{Stating the problem}

Let $4 p = U^2 + D V^2$. We want to find the equation of a curve
$E/\GFq{p}$ having cardinality $m = p+1-U$. The general
algorithm is the following:

\medskip
\noindent
{\bf procedure} {\sc BuildEWithCM($D, U, V, p$)}

\{ Input: $4 p = U^2 + D V^2$ \}

1. For some invariant $u$, compute the minimal polynomial $H_D[u](X)$.

2. Find a root $x_0$ of $H_D[u](X)$ modulo $p$.

3. for all roots $j$ of $\Phi[u](x_0, J)\bmod p$ do

\ \ \ \ a. compute $E(j)$.

\ \ \ \ b. If $\#E(j) = p+1+U$ instead of $p+1-U$, replace $E(j)$ by a twist.

\subsection{Eliminating bad curves}

In general, the degree of $\Phi[u](x_0, J)$ is larger than $1$ and we
expect several roots in $J$, not all of which are invariants of
the curves we are looking for. 

In order to eliminate bad curves, we can use the following
result. First, note that:
$$\Delta(E(j)) = 2^{12}\cdot 3^6 j^2/(j-1728)^3.$$
\begin{prop}\label{prop:g3sq}
Let $4 p = U^2 + D V^2$. The number $\Delta(E(j))$ is a square
modulo $p$ in the following cases:

(i) $D$ odd;

(ii) $4 \mid D$ and $2 \mid V$.
\end{prop}

\noindent
{\em Proof:}

(i) If $\alpha$ is as in Theorem \ref{thm:g3}, we deduce that
$\sqrt{-D}\gamma_3(\alpha)$
is in $\OK$, which means that $H_{-D}[\sqrt{-D}\gamma_3]$
splits modulo $p$ and therefore $j-1728 = -D u^2 \bmod p$ and we have
$\legendre{-D}{p}=+1$ by hypothesis.

(ii) Theorem \ref{thm:g3} tells us that
$\Q(\gamma_3(\alpha)) = \Q(j(2 \alpha))$. But $p$ splits in the order
$\Ot{2}$ and therefore in $\Omega_{2t}$, which shows that the minimal
polynomial of $\gamma_3$ splits modulo $p$, proving the result. $\Box$

Coming back to our problem, we see that when the above result applies,
a good curve is such that $\legendre{\Delta(E(j))}{p}$ must be equal
to $1$.

\subsection{Deciding which curve is good}

We can assume that we are left with only one possible $j$ and
that we want to compute the cardinality of $E(j)$ as quickly as
possible. Let us explain our idea. Let $\mathcal{D} = D V^2$.
Suppose that $\ell\neq p$ is an odd prime (the case $\ell=2$ will be dealt
with later) and $\legendre{-D}{\ell} \neq -1$. In that case,
Theorem \ref{elkies} applies and if we can find one eigenvalue
$\lambda$, we can find $U \bmod \ell$. If $U\not\equiv 0\bmod\ell$,
then we can find the sign of $U$. Note that if $\ell\mid D$,
then $U\not\equiv 0\bmod\ell$.

The most favorable case is when $\ell
\mid\mathcal{D}$, because then there is only one eigenvalue $\lambda$
(it can be a double one) and $\lambda \equiv U/2 \bmod\ell$. Having
$\lambda$ gives us immediately the sign of $U$. A very favorable case
is when $\ell\equiv 3\bmod 4$, using Dewaghe's idea.

Apart from this, there is another interesting sub-case, when
we can find a rational root $x_0$ of $g_{\ell, \lambda}^E$, using for
instance some class invariant. In that case, we can form $y_0^2 = x_0^3+a
x_0+b\bmod p$ and test whether $y_0$ is in $\GFq{p}$ or not. If it is,
then $\lambda=1$, since $(x_0, y_0)$ is rational and $\pi(P)=P$.
Otherwise, $\lambda=-1$.

Our idea is then to use the general framework for some precise values
of $\ell$, and use rational roots of $g_{\ell, \lambda}$ obtained via
class invariants.
When $\ell=3$, we are sure to end with a rational root of
$f_3^{E(j)}(X)$, as is the case for $\ell=2$ and
$f_4^{E(j)}$. Moreover, we can use some invariant that give us the
torsion points directly. We also give examples for $\ell=5, 7, 11$.
\cbend

\section{The case $\ell=3$}

We suppose that $4 p = U^2+D V^2$. The first subsection makes precise the
above results.

\subsection{Using $3$-torsion points}

We begin with an easy lemma that can be proved by algebraic manipulations:
\begin{lemma}\label{threetors}
Let $v$ be any root of $\Phi_3^c(X, j)=0$.
Then a root of $f_3^{E(j)}(X)$ is given by
$$x_3 = -\frac{(v+27) (v+3)}{v^2+18 v-27}.$$
\end{lemma}

\begin{prop}\label{fdalprop}
Let $p$ be a prime representable as $4 p = U_0^2 + D V_0^2$, for which $3
\mid D V_0^2$ and $\#E = p+1-U$.
Suppose $P = (x_3, y_3)$ is a $3$-torsion point on $E(j)$ for which $x_3$
is rational. Let $s = x_3^3+a_4(j) x_3 + a_6(j) \bmod p$.
Then $U \equiv 2 \legendre{s}{p} \bmod 3$.
\end{prop}

\noindent
{\em Proof:} This is a simple application of Theorem \ref{elkies}.
$\Box$

 \subsection{Solving the equation $\Phi_3^c(X, j(E))=0$}

  \subsubsection{The case $\legendre{-D}{3}\neq -1$}

A solution of this equation is given by $\galfone{3}^{12}$, which lies
in $\Omega_1$ with the hypothesis made on $D$. 

\noindent
{\bf Numerical examples.} Let $H_{-15}[\galfone{3}^{12}] = {X}^{2}+81\,X+729$,
$p=109$, $4p = 14^2+15 \times 4^2$, $v_3=3$, $x_3=104$, $E: Y^2=X^3+94
X+99$; $U=\pm 14$. Since $\lambda=1\bmod 3$, we conclude that $U = 14$
and $E$ has $109+1-14$ points.

Take $D = 20$ and $p = 349$. We find $(U, V) = (\pm 26, \pm 6)$. We compute:
$$H_{-20}[\galfone{3}^{12}] = X^2+(70-22 \sqrt{-20}) X-239-154 \sqrt{-20}.$$
Using $\sqrt{-20}=237\bmod p$, a root of this polynomial is $v_3 = 257$, from
which $j=224$ and $E(j): Y^2=X^3+45 X+30$. Now $\lambda=-1$, which
gives us that $\#E = 349+1+26$.

  \subsubsection{The case $\legendre{-D}{3}=-1$}

We may solve the degree $4$ equation $\Phi_3^c(X,
j(\alpha))=0$ directly.

In Skolem's approach \cite{Skolem52a}, to compute the roots of a
general quartic (with $a_1$ and $a_3$ not both zero)
$$P(X)=X^4+a_1 X^3 + a_2 X^2+a_3 X + a_4$$
one uses the four roots $X_i$ of $P$ to define
\begin{equation}\label{Xz}
 \left\{\begin{array}{lcl}
 z_1 &=& X_1 + X_2 - X_3 - X_4, \\
 z_2 &=& X_1 - X_2 + X_3 - X_4, \\
 z_3 &=& X_1 - X_2 - X_3 + X_4.
        \end{array}
 \right.
\end{equation}
Writing $y_i = z_i^2$, the $y_i$'s are roots of
\begin{equation}
  R(y) = y^3 + b_1 y^2+ b_2 y + b_3
\end{equation}
in which
\begin{equation}
 \left\{\begin{array}{ccl}
  b_1 &=& 8 a_2 -3 a_1^2, \\
  b_2 &=& 3 a_1^4 - 16 a_1^2 a_2 + 16 a_1 a_3 + 16 a_2^2 -64 a_4, \\
  b_3 &=& -(a_1^3 - 4 a_1 a_2 + 8 a_3)^2.
 \end{array}\right.
\end{equation}
Conversely, if the $y_i$'s are the roots of $R$ and if the $z_i$'s
are chosen in such a way that
$$- z_1 z_2 z_3 = a_1^3 - 4 a_1 a_2 + 8 a_3,$$
then the $X_i$'s defined by (\ref{Xz}) (together with
$X_1+X_2+X_3+X_4 = -a_1$) are the roots of $P$.

In our case, we find that
$$R(Y) = Y^3 - 1728 Y^2 - 576 (j(\alpha) - 1728) Y -64 (j(\alpha) -
1728)^2$$
and the compatibility relation is $z_1 z_2 z_3 = 8 (j(\alpha) - 1728)$.
Since we suppose that $3 \nmid D$, we replace $j(\alpha)$ by
$\gamma_2(\alpha)^3$. In that case, the roots of $R(Y)$ are
$$4 (\zeta_3^{2i}\gamma_2(\alpha)^2+12 \zeta_3^i\gamma_2(\alpha)+144)$$
for $i = 0, 1, 2$. Studying the roots of these numbers as class invariants
could probably be done. The function
$$\sqrt{\gamma_2(\alpha)^2+12 \gamma_2(\alpha)+144}$$
has been introduced via a different route by Birch in \cite{Birch69}
and the theorems proven there could be used in our context, though we
refrain from doing so in this article.

Let us summarize the algorithm to find the roots of $\Phi_3^c(X,
j(E))$ modulo $p$ when $3\nmid D, 3\mid V$ (which implies $p\equiv
1\bmod 3$):

1. compute $\gamma_2\bmod p$;

2. compute the values $y_i = 4 (\zeta_3^{2i}\gamma_2(\alpha)^2+12
   \zeta_3^{i}\gamma_2(\alpha)+144) \bmod p$ for $i=1, 2$;

3. compute $z_i = \sqrt{y_i}\bmod p$ for $i=1, 2$ and $z_3 = 8
   (\gamma_2^3-1728)/(z_1 z_2)$ from which $X_1 = z_1+z_2+z_3-36$ is a
   root of $\Phi_3^c(X, j)$.

Notice that $\zeta_3\bmod p$ can be computed as follows (see
\cite{Atkin92} for more on this sort of ideas): since $3\mid
p-1$, we can find $a$ such that $a^{(p-1)/3}\not\equiv 1\bmod p$. Put
$\zeta_3 = a^{(p-1)/3}$. It satisfies $\zeta_3^2+\zeta_3+1\equiv
0\bmod p$. Therefore, finding a root costs two squareroots and one
modular exponentiation, once $\gamma_2$ is known.

\noindent
{\bf Numerical examples.} 
Consider $(D, p, U, V)=(40, 139, \pm 14, \pm3)$. A root of
$H_{-40}[\gamma_2](X) = {X}^{2}-780\,X+20880$ modulo $p$ is
$110$. Using $\zeta_3 = 96$, we compute $v_3 = 109$ and $x_3 =
135$. Then $E: Y^2 = X^3 + 124 X+129$ has $\lambda=1$ and $U = 14$.

\section{The case $\ell = 5$}

\subsection{Using $\galfone{5}$}

We assume here that $\legendre{-D}{5} \neq -1$ and $5 \mid D V^2$. In
that case, we can use some power of $\galfone{5}$ as invariant to get
a root $v_5$ of $\Phi_5^c(X, j)$, thus yielding a factor $g_5^{E(j)}$
of $f_5^{E(j)}$. Writing:
$$A = v_5^2+22 v_5+125, B = v_5^2+4 v_5-1, C = v_5^2+10 v_5+5,$$
one has:
$$g_5^{E(j)}(X) = X^2 + 2 (C/B) X + (1-36/A) (C/B)^2.$$
Putting $Y = (B/C) X$ leads us to $(Y+1)^2 - 36/A$. At this point,
since
$$j = \frac{(v_5^{2}+10\,v_5+5)^{3}}{v_5}$$
we also have:
$$j - 1728 = {\frac {\left ({v_5}^{2}+22\,v_5+125\right
)\left ({v_5}^{2}+4\,v_5-1\right )^{2}}{v_5}}$$
or $A = v_5 (j-1728) / B^2$.

\subsubsection{The case $U\equiv \pm 2\bmod 5$}

We deduce that $p\equiv 1\bmod 5$ and $g_5^{E(j)}(X)$ has
two rational roots. 

\noindent
{\bf Examples.} Take $D=35$ for which
$$H_{-35}[\galfone{5}^6](X) = X^2+50 X+125.$$
Take $(p, U, V) = (281, \pm 33, \pm 1)$. We first use $v_5 = 163$ to
compute $E(j): Y^2 = X^3+32 X+115$ and
$g_5^{E(j)}(X) = X^2+245 X+198$. From this, we get $x_5 = 227$ and find that
$x_5^3+a_4(j)x_5+a_6(j)$ is a square in $\GFq{p}$, so that
$\#E(j)=p+1+33$.

Consider now $D=91$ for which $\legendre{-91}{5}=+1$. We find:
$$H_{-91}[\galfone{5}^6] = X^2+(130-40 \sqrt{-91}) X-99-8 \sqrt{-91}.$$
Taking $(p, U, V)=(571, \pm 3, \pm 5)$, we use $\sqrt{-91}=342\bmod p$, find
$v_5 = 216$ from which $j=533$ and $E(j): Y^2 = X^3+181 X+ 311$. Then
$g_5^{E(j)}(X) = X^2+213 X+412$ which has a root $x_5 = 315$. We find
that $\lambda=-1$ and $U = 3$.

\subsubsection{The case $U\equiv \pm 1\bmod 5$}

One has $p\equiv 4\bmod 5$ and $g_5^{E(j)}(X)$ is
irreducible; the eigenvalue is $\lambda=U/2\equiv \pm 2\bmod 5$. We
can compute it using the techniques of SEA, that is test the identity
$$(X^p, Y^p) = [\pm 2](X, Y) \bmod g_5^{E(j)}(X).$$
(Actually, checking the equality on the ordinates is enough.)
Depending on the implementation, this can cost more than testing
$[m]P$ on $E$.

\noindent
{\bf Example.} Consider $(D, p, U, V)=(35, 109, \pm 11, \pm 3)$. One
computes $v_5 = 76$ and $g_5^{E(j)}(X) = X^2+13 X+13$. We compute
$$(X^p, Y^p) \equiv (108 X + 96, Y (72 X+43)) = [2] (X, Y).$$
Therefore, $U = -11$.

Consider $(D, p, U, V) = (91, 569, \pm 1, \pm 5)$. We find $E(j): Y^2
= X^3+558 X + 372$, $g_5^{E(j)}(X) = X^2+100 X+201$ and
$$(X^p, Y^p) \equiv [2] (X, Y)$$
so that $U = -1$.

\subsection{A remark on the case $D=20$}\label{ssct:20}

We will take a route different from that in \cite{LeMo97}.
Write $p = a^2 + 5 b^2$. Let $\varepsilon_0 = (1+\sqrt{5})/2$
be the fundamental unit of $\Q(\sqrt{5})$. We have
$$a_4=-{\frac {162375}{87362}}-{\frac {89505\,\sqrt{5}}{174724}}, \quad
a_6=-{\frac {54125}{43681}}-{\frac {29835\,\sqrt{5} }{87362}}$$
and $f_5(X)$ has the factor:
$${X}^{2}+\left({\frac{695}{418}}+{\frac {225\,\sqrt{5}}{418}}\right )X+
{\frac {129925}{174724}}+{\frac{45369\,\sqrt{5}}{87362}}
$$
of discriminant:
$$\Delta = \frac{3^2}{11^2 \cdot 19^2}
\left(\frac{7+\sqrt{5}}{2}\right)^4
\left(\frac{9+\sqrt{5}}{2}\right)^2 \frac{\sqrt{5}}{\varepsilon_0^5}$$
which is congruent to $\varepsilon_0 \sqrt{5}$ modulo squares. Now, by
\cite{Lehmer72}, we have
$$\myjac{\varepsilon_0 \sqrt{5}}{p} = \myjac{p}{5}_4.$$
When $p \equiv 1 \bmod 20$, $\Delta$ is a square modulo $p$ and there
are two abscissas in $\GFq{p}$. Now, $a \equiv \pm 1 \bmod 5$ and
thus
$$\#E(j) \equiv 1 + 1 \pm 2 \bmod 5.$$
We can distinguish the two cases by computing $y_5$: It is in
$\GFq{p}$ if and only if $m \equiv 0\bmod 5$.

\section{Numerical examples for $\ell \equiv 3\bmod 4$}

\subsection{The case $\ell = 7$}

\begin{lemma}
Let $v_7$ be a root of $\Phi_7^c(X, j)$ and put
$$A(v_7) = v_7^{4}+14\,v_7^{3}+63\,v_7^{2}+70\,v_7-7.$$
Then
$$\mathrm{Resultant}(g_{7, \lambda}^{E(j)}(X), X^3+a_4(j)X+a_6(j)) = -3
j v_7 A(v_7) S(v_7)^2$$
for some rational fraction $S$ with integer coefficients.
\end{lemma}

\noindent
{\em Proof:} using {\sc Maple}, we compute:
$$\mathrm{Resultant}(g_{7, \lambda}, X^3+a_4(j)X+a_6(j)) 
= -2^{12} \cdot 3^9 \cdot \left ({v_7}^{2}+13\,v_7+49\right
)^{3}\left ({v_7}^{2}+5\,v_7+1\right )^{9}/A^9$$
from which the result follows. $\Box$

Take $D=91$ for which
$$H_{-91}[\galfone{7}^4] = X^2+77 X+49.$$
Take $(p, U, V)=(107, \pm 8, \pm 2)$. We find $v_7 = 62$ from which
$g_7^{E(j)}(X) = X^3+104 X^2+44 X+73$. Using $E(j): Y^2=X^3+101
X+103$, we find $r = 13$ and $\legendre{13}{p}=1$ and therefore $U =
8$.

For $(D, p, U, V) = (20, 569, \pm 36, \pm 7)$, we compute:
$$H_{-20}[\galfone{7}^4](X)=X^2+(15-\sqrt{-20}) X+41-6 \sqrt{-20}$$
one of which roots modulo $p$ is $v_7 = 195$ (taking
$\sqrt{-20}=320$). Then $E(j): Y^2 = X^3+289 X+ 3$ has $g_7^{E(j)}(X)
= X^3+111 X^2+185 X+94$ from which $U = 36$.

\subsection{The case $\ell = 11$}

In that case, the modular equation is quite large. However, if we
restrict to the case where $3\nmid D$, we can use the modular equation
relating $\galfone{11}^4$ and $\gamma_2$:
{
\small
\begin {eqnarray*}
&&
{X}^{12}-1980\,{X}^{9}+880\,\gamma_2{X}^{8}+44\,\gamma_2^{2}{X}^{7}+98007
8\,{X}^{6}-871200\,\gamma_2{X}^{5}+150040\,\gamma_2^{2}{X}^{4} \\
&& +\left (47066580-7865\,\gamma_2^{3}\right ){X}^{3}
+\left(154\,\gamma_2^{4}+560560\,\gamma_2\right ){X}^{2}+\left
(1244\,\gamma_2^{2}-\gamma_2^{5}\right )X+121.
\end {eqnarray*}
}
Consider $(D, p, U, V) = (88, 103, \pm 18, \pm 1)$. First, we find:
$$H_{-88}[\galfone{11}^4](X)=X^2-66 X+121$$
a root of which is $w_{11} = 21$. Plugging this into the modular
equation, we find $\gamma_2 = 63$, from which $j = 66$ and $E(j): Y^2
= X^3 + 73 X+ 83$. Using the techniques of SEA, we find that
$$g_{11} = X^5+81 X^4+22 X^3+55 X^2+99 X+15$$
and the resultant is $98$, so that $U=18$.

Note that the techniques needed to compute $g_{11}$ are probably too
heavy to make this case useful. However, we provide it as a non-trivial
example.

\section{The case $\ell = 2$}

The points of $2$-torsion cannot be used in our context, since they
have ordinate $0$. So we must try to use $4$-torsion points
instead. We suppose that $-D$ is fundamental.

\subsection{Splitting $f_4^{E(j)}$}

Curves having rational $2$-torsion are parametrized by $X_0(2)$, or
equivalently, $j(E) = (u+16)^3/u$. Notice that:
\begin{equation}\label{j1728}
j-1728 = \gamma_3^2 = {\frac {\left (u+64\right )\left (u-8\right
)^{2}}{u}}.
\end{equation}
Using algebraic manipulations (and {\sc
Maple}), $f_4^{E(j)}(X)$ factors as $P_2(X) P_4(X)$ where:
$$P_2(X) = X^2 + 2 \,{\frac {u+16}{u-8}} X + {\frac {\left
(u-80\right )\left (u+16\right )^{2}}{\left (u-8\right )^{2}\left
(u+64\right )}},$$
$$P_4(X) = {X}^{4}
-2\,{\frac {u+16}{u-8}} {X}^{3}
-12\,{\frac {\left (u+16\right )^{2}}{\left (u+64\right )\left
(u-8\right )}}{X}^{2}
-2\,{\frac {\left (7\,u+16\right )\left (u+16\right )^{3}}{\left (u+
64\right )\left (u-8\right )^{3}}} X$$
$$-{\frac {\left (5\,{u}^{2}+640\,u-
256\right )\left (u+16\right )^{4}}{\left (u+64\right )^{2}\left (u-8
\right )^{4}}}.$$
The polynomial $P_2$ has discriminant:
$$\Delta_2(u) = 12^2\,{\frac {\left (u+16\right )^{2}}{\left (u-8\right
)^{2}\left (u+64\right )}}.$$
The polynomial $P_4$ has the following property. If $(u+64)/u = v^2$,
then it splits as a product of two quadratic polynomials:
$$G_a(X) = {X}^{2}+2\,{\frac {\left ({v}^{2}+3\right )}{v\left(v+3\right )}}X
+{\frac {\left ({v}^{2}+12\,v-9\right )\left ({v}^{2}+3\right )^{2}}{
\left (v+3\right )^{2}\left (v-3\right )^{2}{v}^{2}}},$$
$$G_b(X) = {X}^{2}+2\,{\frac {\left ({v}^{2}+3\right )}{v\left(v-3\right )}}X
+{\frac {\left ({v}^{2}-12\,v-9\right )\left ({v}^{2}+3\right )^{2}}{
\left (v+3\right )^{2}\left (v-3\right )^{2}{v}^{2}}}.$$

\begin{prop}\label{prop2}
Suppose that $(D, p, V)$ satisfies one of the conditions of Proposition
\ref{prop:g3sq} and that $u$ is a square. Then $P_2$ splits modulo $p$.
\end{prop}

\noindent
{\em Proof:} Equation (\ref{j1728}) tells us that $u (u+64)$ is a
square modulo $p$, which implies that $\Delta_2(u)$
is also a square. $\Box$

Notice that generally, at least one of the roots of $\Phi_2^c(X, j)$,
denoted by $u$, will be the square of some Weber function, see
\cite{Schertz02}.

\subsection{Eigenvalues modulo $2^k$}

Our idea is to use the roots of the characteristic polynomial
$\chi(X) = X^2 - U X + p$ modulo powers of $2$ and deduce from this
the sign of $U$ when possible. This subsection is devoted to
properties of these roots.

Since $p \equiv 1\bmod 2$, $\chi(X)$ has roots modulo $2$ if and only if $U
\equiv 0\bmod 2$. Modulo $4$, $\chi(X)$ has roots if and only if $U
\equiv (p+1) \bmod 4$, which we suppose from now on. It is not enough
to look at this case, since we have $U \equiv 0\bmod 4$ or $U \equiv 2
\bmod 4$ and in both cases, and we cannot deduce from this alone the
sign of $U$. We will need to look at what happens modulo $8$. We list
below the cases where $\chi(X)$ has roots modulo 8 and then relate
this with the splitting of $p$.

\begin{lemma}\label{lem1}
The solutions of $X^2 \equiv 4\bmod 8$ are $\pm 2$.
\end{lemma}

\begin{lemma}\label{lem2}
Write $\varepsilon = \pm 1$. We give in the following table the roots
of $\chi(X)$ modulo $8$:
$$\begin{array}{|c|c|c|c|}\hline
p \bmod 8 \backslash U \bmod 8 & 0 & 2 \varepsilon & 4 \\ \hline
1 & \emptyset & \{\varepsilon, \varepsilon+4\} & \emptyset\\ \hline
3 & \emptyset & \emptyset & \{\pm 1, \pm 3\} \\ \hline
5 & \emptyset & \{-\varepsilon, -\varepsilon+4\} & \emptyset\\ \hline
7 & \{\pm 1, \pm 3\} & \emptyset & \emptyset\\ \hline
\end{array}$$
\end{lemma}

\begin{prop}
Let $4 p = U^2 + D V^2$. The polynomial $\chi(X)$ has roots modulo $8$
exactly in the following cases:

(i) $4 \mid D$ and $2 \mid V$;

(ii) $4 \nmid D$ and [($4 \mid V$) or ($2 \mid\mid V$ and $D \equiv
7\bmod 8$)].
\end{prop}

\noindent
{\em Proof:} 

(i) If $V$ is even, we deduce that $U^2 \equiv 4 p \equiv 4 \bmod
8$, $\chi(X)$ is one of
$X^2-2\varepsilon X+1$ or $X^2-2\varepsilon X+5$ by Lemma
\ref{lem1}. The result follows from Lemma \ref{lem2}.

What can be said when $V$ is odd?
When $4 \mid\mid D$, this means that $p = (U/2)^2+(D/4)V^2$, implying
that $U \equiv 0 \bmod 4$ and $p \equiv 1\bmod 4$ (since $-D$ is
fundamental, $D/4\equiv 1\bmod 4$), but then $U \not\equiv p+1\bmod 4$.

When $8 \mid D$, then $p = (U/2)^2 +(D/4) V^2$ with $U \equiv \pm 2
\bmod 8$, but $p \equiv 3 \bmod 4$ and again $U \not\equiv
p+1\bmod 4$.

(ii) In that case, $U$ and $V$ have the
same parity. If $U$ and $V$ are odd, this implies
$m=p+1-U$ is odd, so that we do not have $2$-torsion points. If $U$
and $V$ are even, so is $m$ and $p = (U/2)^2 + D (V/2)^2$. 

If $V/2$ is even of the form $2 V'$, then $p = (U/2)^2+4 D {V'}^2$;
$U/2$ must be odd and $p \equiv 1\bmod 4$ and we conclude as in case (i).

If $V/2$ is odd, then $p = (U/2)^2+D {V'}^2$ with $V'$ odd, which
implies $U/2$ even, that is $U \equiv 0\bmod 8$ or $U \equiv 4 \bmod
8$. One has $p \equiv (U/2)^2+D\bmod 8$. If $D\equiv 7\bmod 8$, then
$(U, p) = (0, 7) \bmod 8$ or $(4, 3)\bmod 8$ and the two
characteristic polynomials have four roots modulo $8$.
If $D \equiv 3\bmod 8$, then $(U, p) = (0, 3)$ or $(4, 7)$ modulo $8$
and $\chi(X)$ has no roots. $\Box$

\subsection{Computing the cardinality of CM-curves}

This section makes use of the theory of
isogeny cycles described in \cite{CoMo94,CoDeMo96}. 

With the notations of the preceding section,
we suppose we are in the case where $U = 2\varepsilon \bmod 8$, or
equivalently $4\mid D$ and $2\mid V$, or $4\nmid D$ and $4 \mid
V$.

From Proposition \ref{prop2}, we know that the factor $P_2(X)$ of
$f_4^{E(j)}$ has at least two roots modulo $p$. If $x_4$ is one of
these and $s = x_4^3+a x_4+b$, we let $y_4 = \sqrt{s}$ ({\it a priori} in
$\GFpn{p}{2}$) and $P=(x_4,
y_4)$. Now $\pi(P) = \pm P$ according to the fact that $s$ is a square
or not. We have our eigenvalue $\lambda_4 \equiv \pm 1 \bmod 4$. By the
theory of isogeny cycles, the eigenspace $C_4$ generated by $P$ can be
lifted to an eigenspace $C_8$ of $E[8]$ associated to the eigenvalue
$\lambda_8$ which is congruent to $\lambda_4$ modulo $4$. Since $U =
2\varepsilon \bmod 8$, we know from Lemma \ref{lem2} that only one of
the possible values of $\lambda_8$ reduces to a given $\lambda_4$,
which gives us $\varepsilon$.

In practice, $x_4$ is relatively inexpensive to use when $u$ is the
square of a Weber function, which happens in the case $4\mid D$ or $D
\equiv 7\bmod 8$ (for this, one uses an invariant for $-4 D$ instead
of $-D$, and both class groups have the same class number, see
\cite{AtMo93b}). When $D \equiv 3\bmod 4$, $h_t = 3 h_1$, which is not
as convenient; still, a root of $\Phi_2^c(X, j)$ exists, since it is
in $\Omega_2$ and $p$ splits in it.

\noindent
{\bf Examples.} First take $(D, p, U, V)=(20, 29, \pm 6, \pm 2)$. We
find $u=7$, $j=23$ and $E(j): Y^2 = X^3+3 X+2$. From this, $P_2$ has a
root $x_4 = 7$ and $\lambda_8 = -1$, so that $U = -6$.

Now take $(D, p, U, V)=(40, 41, \pm 2, \pm 2)$. We compute $u = 16$,
$j = 39$, $E(j): Y^2 = X^3+30 X+20$, $x_4 = 19$ and $\lambda_8 = -1$
implying $U = -2$.

Let us turn to odd $D$'s. Take $(D, p, U, V) = (15, 409, \pm 26, \pm
8)$. Then $u=102$, $j=93$, $E: Y^2 = X^3+130 X+223$, $x_4 = 159$
yielding $\lambda_8 = -1$ and $U = -26$.

\subsection{The case $D$ odd}

In that case, $\Phi_2^c(X, J)$ will
have three roots in $\GFq{p}$ or $\GFpn{p}{2}$, that we can compute
directly. This could be useful for the cases not treated by the the
preceding section.

Let us try to solve the equation
$$\Phi_2^c(X, J) = {X}^{3}+48\,{X}^{2}+768\,X-J X+4096 = 0$$
directly. As in \cite{Cailler08} (already used in \cite{Morain90c}),
we first complete the cube letting $Y = X + 16$ to get:
\begin{equation}\label{eqY}
{Y}^{3}-JY+16\,J = 0.
\end{equation}
We look for $\alpha$ and $\beta$ such that this equation can be
rewritten:
$$Y^3-3 \alpha \beta Y + \alpha \beta (\alpha + \beta) \equiv 0.$$
The coefficients $\alpha$ and $\beta$ are solutions of
$$W^2 - 48 W + J/3 = 0$$
whose discriminant is $\Delta = (-4/3) (J-1728)$. Having $\alpha$ and
$\beta$ (in $\GFq{p}$ or $\GFpn{p}{2}$), we solve
$$z^3 = \frac{\alpha}{\beta}$$
and we get a root
$$Y = \frac{\beta z - \alpha}{z - 1}$$
of (\ref{eqY}).

\medskip
Since $D$ is odd,
$\sqrt{-D}\gamma_3$ is an invariant, so that we can write:
$$\Delta = -\frac{4}{3}
\left(\frac{\sqrt{-D}\gamma_3}{\sqrt{-D}}\right)^2.$$
The computation of the roots then depends on $\legendre{-3}{p}=1$.
It is not clear that the above mentionned approach is really faster
than the naive one.

\section{Applications to ECPP and conclusion}

In ECPP, the situation is as follows. We are given $j$ and $m=p+1-U$
for some known $U$. We have to build an elliptic curve $E$ having
invariant $j$ and cardinality $m$. We use the results of the preceding
sections in the following way. We build a candidate $E$ and compute
its cardinality $m'$. If $m' = m$, then $E$ is the correct answer,
otherwise, we have to twist it.

In \cite{EnMo02}, a comparison of all possible class invariants for a
given $D$ was made using the height of their minimal
polynomial. Though it is clear that it is easier to use
invariants of small height, the results of the present article show
that we might as well favor those invariants that give us a fast way
of computing the right equation instead.

For instance, if $(D, 6)=1$, using Stark's ideas whenever possible is
a good thing. When $3\mid D$ or $7\mid D$, $\galfone{3}$ or
$\galfone{7}$ should be preferred since we have a fast
answer. Note now a new phenomenon. If we are interested in a
prescribed $p$, we should use an invariant which depends on $D$, but
also on $p$, or more precisely on the small factors of $V$. For
instance, if $3 \mid V$, we can use the direct solution of
$\Phi_3^c(X, J)$. If not, we may use some case where
$\legendre{-D}{\ell}=+1$, and $\ell\mid V$.

The present work has enlarged the set of $D$'s for which the
corresponding $E$'s are easy to find. Nevertheless, there are cases
which are badly covered (for instance odd primes which are non
quadratic residues modulo $8$, $3$, $5$, $7$, such as $D=163$) and
that will require new ideas to be treated.

\medskip
\noindent
{\bf Acknowledgments.} The author wants to thank A.~Enge for his
careful reading of the manuscript and suggesting many improvements.

\iffalse
\begin{small}
\bibliographystyle{plain}
\bibliography{morain}

\begin{thebibliography}{10}

\bibitem{Atkin92b}
A.~O.~L. Atkin.
\newblock The number of points on an elliptic curve modulo a prime ({II}).
\newblock \Draft. Available on {\tt
  http://listserv.nodak.edu/archives/nmbrthry.html}, 1992.

\bibitem{Atkin92}
A.~O.~L. Atkin.
\newblock Probabilistic primality testing.
\newblock In P.~Flajolet and P.~Zimmermann, editors, {\em Analysis of
  Algorithms Seminar {I}}. INRIA Research Report XXX, 1992.
\newblock Summary by F. Morain. Available as {\tt
  http://pauillac.inria.fr/algo/seminars/sem91-92/atkin.ps}.

\bibitem{AtMo93b}
A.~O.~L. Atkin and F.~Morain.
\newblock Elliptic curves and primality proving.
\newblock {\em Math. Comp.}, 61(203):29--68, July 1993.

\bibitem{Birch69}
B.~J. Birch.
\newblock Weber's class invariants.
\newblock {\em Mathematika}, 16:283--294, 1969.

\bibitem{Cailler08}
C.~Cailler.
\newblock Sur les congruences du troisi\`eme degr\'e.
\newblock {\em Enseign. Math.}, 10:474--487, 1902.

\bibitem{CoDeMo96}
J.-M. Couveignes, L.~Dewaghe, and F.~Morain.
\newblock Isogeny cycles and the {S}choof-{E}lkies-{A}tkin algorithm.
\newblock Research Report LIX/RR/96/03, LIX, April 1996.
\newblock Available at {\tt
  http://www.lix.polytechnique.fr/Labo/Francois.Morain/}.

\bibitem{CoMo94}
J.-M. Couveignes and F.~Morain.
\newblock Schoof's algorithm and isogeny cycles.
\newblock In L.~Adleman and M.-D. Huang, editors, {\em Algorithmic Number
  Theory}, volume 877 of {\em Lecture Notes in Comput. Sci.}, pages 43--58.
  Springer-Verlag, 1994.
\newblock 1st Algorithmic Number Theory Symposium - Cornell University, May
  6-9, 1994.

\bibitem{Cox89}
D.~A. Cox.
\newblock {\em Primes of the form $x^2+n y^2$}.
\newblock John Wiley \& Sons, 1989.

\bibitem{Dewaghe98}
L.~Dewaghe.
\newblock Remarks on the {S}choof-{E}lkies-{A}tkin algorithm.
\newblock {\em Math. Comp.}, 67(223):1247--1252, July 1998.

\bibitem{Elkies98}
N.~D. Elkies.
\newblock Elliptic and modular curves over finite fields and related
  computational issues.
\newblock In D.~A. Buell and J.~T. Teitelbaum, editors, {\em Computational
  Perspectives on Number Theory: Proceedings of a Conference in Honor of A. O.
  L. Atkin}, volume~7 of {\em AMS/IP Studies in Advanced Mathematics}, pages
  21--76. American Mathematical Society, International Press, 1998.

\bibitem{EnMo02}
A.~Enge and F.~Morain.
\newblock Comparing invariants for class fields of imaginary quadratic fields.
\newblock In C.~Fieker and D.~R. Kohel, editors, {\em Algorithmic Number
  Theory}, volume 2369 of {\em Lecture Notes in Comput. Sci.}, pages 252--266.
  Springer-Verlag, 2002.
\newblock 5th International Symposium, ANTS-V, Sydney, Australia, July 2002,
  Proceedings.

\bibitem{Ishii04}
N.~Ishii.
\newblock Trace of {F}robenius endomorphism of an elliptic curve with complex
  multiplication.
\newblock Available at {\tt http://arxiv.org/abs/math.NT/0401289}, January
  2004.

\bibitem{JoMo95}
A.~Joux and F.~Morain.
\newblock Sur les sommes de caract{\`e}res li{\'e}es aux courbes elliptiques
  {\`a} multiplication complexe.
\newblock {\em J. Number Theory}, 55(1):108--128, November 1995.

\bibitem{Kohel96}
D.~Kohel.
\newblock {\em Endomorphism rings of elliptic curves over finite fields}.
\newblock PhD thesis, University of California at Berkeley, 1996.

\bibitem{Lehmer72}
E.~Lehmer.
\newblock On some special quartic reciprocity law.
\newblock {\em Acta Arith.}, {XXI}:367--377, 1972.

\bibitem{LeMo97}
F.~Lepr{\'e}vost and F.~Morain.
\newblock Rev\^{e}tements de courbes elliptiques \`{a} multiplication complexe
  par des courbes hyperelliptiques et sommes de caract\`{e}res.
\newblock {\em J. Number Theory}, 64:165--182, 1997.

\bibitem{MaMu01}
M.~Maurer and V.~M{\"u}ller.
\newblock Finding the eigenvalue in {E}lkies' algorithm.
\newblock {\em Experiment. Math.}, 10(2):275--285, 2001.

\bibitem{Morain90c}
F.~Morain.
\newblock {\em Courbes elliptiques et tests de primalit{\'e}}.
\newblock Th{\`e}se, Universit{\'e} Claude Bernard--Lyon {I}, September 1990.

\bibitem{Morain95a}
F.~Morain.
\newblock Calcul du nombre de points sur une courbe elliptique dans un corps
  fini~: aspects algorithmiques.
\newblock {\em J. Th\'eor. Nombres Bordeaux}, 7:255--282, 1995.

\bibitem{Morain98a}
F.~Morain.
\newblock Primality proving using elliptic curves: an update.
\newblock In J.~P. Buhler, editor, {\em Algorithmic Number Theory}, volume 1423
  of {\em Lecture Notes in Comput. Sci.}, pages 111--127. Springer-Verlag,
  1998.
\newblock Third International Symposium, ANTS-III, Portland, Oregon, june 1998,
  Proceedings.

\bibitem{Newman57}
M.~Newman.
\newblock Construction and application of a class of modular functions.
\newblock {\em Proc. London Math. Soc.}, 3(7):334--350, 1957.

\bibitem{Newman59}
M.~Newman.
\newblock Construction and application of a class of modular functions ({II}).
\newblock {\em Proc. London Math. Soc.}, 3(9):373--387, 1959.

\bibitem{PaVe96}
R.~Padma and S.~Venkataraman.
\newblock Elliptic curves with complex multiplication and a character sum.
\newblock {\em J. Number Theory}, 61:274--282, 1996.

\bibitem{Schertz02}
R.~Schertz.
\newblock Weber's class invariants revisited.
\newblock {\em J. Th\'eor. Nombres Bordeaux}, 14:325--343, 2002.

\bibitem{Schoof95}
R.~Schoof.
\newblock Counting points on elliptic curves over finite fields.
\newblock {\em J. Th\'eor. Nombres Bordeaux}, 7:219--254, 1995.

\bibitem{Silverman86}
J.~H. Silverman.
\newblock {\em The arithmetic of elliptic curves}, volume 106 of {\em Grad.
  Texts in Math.}
\newblock Springer, 1986.

\bibitem{Silverman94}
J.~H. Silverman.
\newblock {\em Advanced Topics in the Arithmetic of Elliptic Curves}, volume
  151 of {\em Grad. Texts in Math.}
\newblock Springer-Verlag, 1994.

\bibitem{Skolem52a}
Th. Skolem.
\newblock The general congruence of 4th degree modulo $p$, $p$ prime.
\newblock {\em Norsk. Mat. Tidsskr}, 34:73--80, 1952.

\bibitem{Stark96}
H.~M. Stark.
\newblock Counting points on {CM} elliptic curves.
\newblock {\em Rocky Mountain J. Math.}, 26(3):1115--1138, 1996.

\bibitem{Swan62}
R.~G. Swan.
\newblock Factorization of polynomials over finite fields.
\newblock {\em Pacific J. Math.}, 12:1099--1106, 1962.

\bibitem{Velu71}
Jacques V{\'e}lu.
\newblock Isog{\'e}nies entre courbes elliptiques.
\newblock {\em C. R. Acad. Sci. Paris S\'er. I Math.}, 273:238--241, July 1971.
\newblock S{\'e}rie A.

\bibitem{Weber02}
H.~Weber.
\newblock {\em Lehrbuch der Algebra}, volume {I}, {II}, {III}.
\newblock Chelsea Publishing Company, New York, 1902.

\end{thebibliography}
\end{small}
\else
\def\noopsort#1{}\ifx\bibfrench\undefined\def\biling#1#2{#1}\else\def\biling#1%
#2{#2}\fi\def\Inpreparation{\biling{In preparation}{en
  pr{\'e}paration}}\def\Preprint{\biling{Preprint}{pr{\'e}version}}\def\Draft{%
\biling{Draft}{Manuscrit}}\def\Toappear{\biling{To appear}{\`A para\^\i
  tre}}\def\Inpress{\biling{In press}{Sous presse}}\def\Seealso{\biling{See
  also}{Voir {\'e}galement}}\def\Editor{\biling{Ed.}{R{\'e}d.}}

\fi
\end{document}